\documentclass[11pt, twoside]{article}
\usepackage{latexsym}
\usepackage{amsmath}
\usepackage{amssymb}
\usepackage[all]{xy}
\usepackage{amsfonts}
\usepackage{verbatim}
\usepackage{amsthm}
\usepackage{mathrsfs}
\usepackage{epsfig}
\usepackage{xy}
\usepackage{array}
\usepackage{stmaryrd}
\usepackage{graphicx,color}
\usepackage{xcolor}
\usepackage[colorlinks=true,linkcolor=blue,citecolor=blue]{hyperref}
\usepackage{tikz}
\usetikzlibrary{arrows,calc}
\usepackage{etex}
\usepackage{mathdots}
\usepackage{float}
\usepackage{graphics}
\usepackage{pdflscape}
\usepackage{bm}
\usepackage{anysize,hyperref}
\input xypic
\xyoption{all}

\usepackage[perpage,symbol]{footmisc}
\usepackage{setspace}
\def\A{\mathscr{A}}

\def\P{\mathcal{P}}
\def\I{\mathcal{I}}

\def\C{\mathscr{C}}
\def\E{\mathbb{E}}
\def\F{\mathbb{F}}
\def\s{\mathfrak{s}}\def\t{\mathfrak{t}}
\def\id{\mathrm{id}}
\def\op{^\mathrm{op}}
\def\Ab{\mathit{Ab}}
\def\del{\delta}
\def\dr{\ar@{->}[r]}

\def\X{\mathscr{X}}\def\T{\mathscr{T}}

\def\Hom{\mbox{Hom}}

\newcommand{\CC}{{\bf{C}}^{n+2}_{\C}}

\newcommand{\ov}{\overset}
\newcommand{\lra}{\longrightarrow}
\newcommand{\co}{\colon}
\newcommand{\uas}{^{\ast}}            
\newcommand{\sas}{_{\ast}}
\newcommand{\Xd}{\langle X_{\bullet},\del\rangle}  
\newcommand{\Yr}{\langle Y_{\bullet},\rho\rangle}  

\newcommand{\ush}{^\sharp}           
\newcommand{\ssh}{_\sharp}
\SelectTips{cm}{10}

\begin{document}
\title{\Large{\bf Hereditary $\bm{n}$-exangulated categories\footnotetext{\hspace{-1em} ~Jian He was supported by the National Natural Science Foundation of China (Grant No. 12171230) and Youth Science and Technology Foundation of Gansu Provincial (Grant No. 23JRRA825).
Jing He was supported by the Hunan Provincial Natural Science Foundation of China (Grant No. 2023JJ40217). Panyue Zhou was supported by the National Natural Science Foundation of China (Grant No. 12371034) and the Hunan Provincial Natural Science Foundation of China (Grant No. 2023JJ30008).}}}
\medskip
\author{Jian He, Jing He and Panyue Zhou}

\date{}

\maketitle
\def\blue{\color{blue}}
\def\red{\color{red}}

\newtheorem{theorem}{Theorem}[section]
\newtheorem{lemma}[theorem]{Lemma}
\newtheorem{corollary}[theorem]{Corollary}
\newtheorem{proposition}[theorem]{Proposition}
\newtheorem{conjecture}{Conjecture}
\theoremstyle{definition}
\newtheorem{definition}[theorem]{Definition}
\newtheorem{question}[theorem]{Question}
\newtheorem{remark}[theorem]{Remark}
\newtheorem{remark*}[]{Remark}
\newtheorem{example}[theorem]{Example}
\newtheorem{example*}[]{Example}
\newtheorem{condition}[theorem]{Condition}
\newtheorem{condition*}[]{Condition}
\newtheorem{construction}[theorem]{Construction}
\newtheorem{construction*}[]{Construction}

\newtheorem{assumption}[theorem]{Assumption}
\newtheorem{assumption*}[]{Assumption}

\baselineskip=17pt
\parindent=0.5cm
\vspace{-6mm}

\begin{abstract}
\baselineskip=16pt
Herschend--Liu--Nakaoka introduced the concept of $n$-exangulated categories as higher-dimensional analogues of extriangulated categories defined by Nakaoka--Palu. The class of $n$-exangulated categories contains $n$-exact categories and $(n+2)$-angulated categories as specific examples. In this article, we introduce the notion of hereditary $n$-exangulated categories, which generalize hereditary extriangulated categories. We provide two classes of hereditary $n$-exangulated categories through closed subfunctors. Additionally, we define the concept of $0$-Auslander $n$-exangulated categories and discuss the circumstances under which these two classes of hereditary $n$-exangulated categories become $0$-Auslander.\\[0.2cm]
\end{abstract}

\pagestyle{myheadings}
\markboth{\rightline {\scriptsize   Jian He, Jing He and Panyue Zhou}}
         {\leftline{\scriptsize Hereditary $n$-exangulated categories}}

\section{Introduction}
The notion of extriangulated categories was introduced by Nakaoka and Palu in \cite{NP} as a simultaneous generalization of
exact categories and triangulated categories. Exact categories (abelian categories are also exact categories) and triangulated categories are extriangulated categories, while there exist some other examples of extriangulated categories which are neither exact nor triangulated, see \cite{NP,ZZ,HZZ,ZhZ,NP1,FHZZ}.

In \cite{GKO}, Geiss, Keller, and Oppermann introduced a new category known as $(n+2)$-angulated categories, which serve as generalizations of triangulated categories; the classical triangulated categories represent the special case when $n=1$. These categories, for instance, emerge as $n$-cluster tilting subcategories of triangulated categories that are closed under the $n$th power of the shift functor.
Subsequently, Jasso \cite{Ja} introduced $n$-exact categories, characterized by certain exact sequences with $n+2$ terms, termed $n$-exact sequences. The case $n=1$ aligns with the usual concepts of exact categories.
An important source of examples for $n$-exact categories is provided by $n$-cluster tilting subcategories of exact categories, as illustrated in \cite[Theorem 4.14]{Ja}.

Recently, Herschend, Liu, and Nakaoka \cite{HLN} introduced the concept of $n$-exangulated categories. It is important to note that the case $n=1$ corresponds to extriangulated categories. Notably, $n$-exact categories ($n$-abelian categories are also $n$-exact categories) and $(n+2)$-angulated categories serve as typical examples of $n$-exangulated categories, as demonstrated in \cite[Proposition 4.5 and Proposition 4.34]{HLN}. Additionally, there exist $n$-exangulated categories that do not fall into the categories of $n$-exact or $(n+2)$-angulated, as discussed in \cite{HLN, HLN1, LZ, HZZ1}.

Gorsky, Nakaoka, and Palu \cite{GNP} recently introduced the concepts of hereditary extriangulated categories and $0$-Auslander extriangulated categories. They developed a theory of mutation for silting objects in certain hereditary extriangulated categories. Additionally, they demonstrated that various notions of maximality for rigid subcategories become equivalent in $0$-Auslander extriangulated categories.
Naturally, we extend these notions to define hereditary $n$-exangulated categories and $0$-Auslander $n$-exangulated categories, as introduced in Definition \ref{def1} and Definition \ref{def2}.

Our first main result is the following.

\begin{theorem}{\rm (see Theorem \ref{main1} for details)}
Let $(\C,\E,\s)$ be an $n$-exangulated category and
$\X$ be a cluster-tilting subcategory of $\C$. We define $$\hspace{4mm}\E_{\X}(C,A)=\{\delta\in\E(C,A)\mid(\del\ssh)_X=0~\mbox{for any}~X\in\X\},$$
$$\E^{\X}(C,A)=\{\delta\in\E(C,A)\mid\del\ush_X=0~\mbox{for any}~X\in\X\}.$$
Then $(\C,\E_{\X},\s\hspace{-1.2mm}\mid_{\E_{\X}})$ and $(\C,\E^{\X},\s\hspace{-1.2mm}\mid_{\E^{\X}})$ are two hereditary $n$-exangulated categories.
\end{theorem}

Let $(\C,\Sigma^n,\Theta)$ be an $(n+2)$-angulated category. It can be regarded as an $n$-exangulated category, denoted by $(\C,\E_{\Sigma^n},\s_{\Theta})$. Assume that $\T$ is a cluster-tilting subcategory of $\C$. We define
$$\F_{\T}(C,A)=\{A\xrightarrow{~~}X_1\xrightarrow{~~}\cdots\xrightarrow{~~}
X_{n-1}\xrightarrow{~~}X_n\xrightarrow{~~}C\overset{\delta}{\dashrightarrow}~\mid~\del~ \text{is induced by the}~(n+2)\text{-angle} $$
\vspace{-8mm}
$$A\xrightarrow{~~}X_1\xrightarrow{~~}\cdots\xrightarrow{~~}X_{n-1}\xrightarrow{~~}X_n\xrightarrow{~~}C\xrightarrow{~f~}\Sigma^nA~\text{with}~f~\text{factors through}~\Sigma^nT,\hspace{1mm} T\in\T\}. $$

Our second main result is the following.

\begin{theorem}{\rm (see Theorem \ref{main2} for details)}
Let $(\C,\Sigma^n,\Theta)$ be an $(n+2)$-angulated category and
$\T$ be a cluster-tilting subcategory of $\C$.
Then

{\rm (1)}~ $(\C,\F_{\T},\s_{\Theta}\hspace{-1.2mm}\mid_{\F_{\T}})$
is a $n$-exangulated category.

{\rm (2)}~ If $\P\subseteq\T$, where $\P$ denotes the full subcategory of projective objects in $(\C,\F_{\T},\s_{\Theta}\hspace{-1.2mm}\mid_{\F_{\T}})$, then $(\C,\F_{\T},\s_{\Theta}\hspace{-1.2mm}\mid_{\F_{\T}})$
is $0$-Auslander.
\end{theorem}

This article is organized as follows. In Section 2, we review elementary definitions and facts on $n$-exangulated categories. In Section 3, we present the proof of our first main result. In Section 4, we give the proof of our second main result.

\section{Preliminaries}
In this section, let $\C$ be an additive category and $n$ be a positive integer. Suppose that $\C$ is equipped with an additive bifunctor $\E\colon\C\op\times\C\to{\rm Ab}$, where ${\rm Ab}$ is the category of abelian groups. Next we briefly recall some definitions and basic properties of $n$-exangulated categories from \cite{HLN}. We omit some
details here, but the reader can find them in \cite{HLN}.

{ For any pair of objects $A,C\in\C$, an element $\del\in\E(C,A)$ is called an {\it $\E$-extension} or simply an {\it extension}. We also write such $\del$ as ${}_A\del_C$ when we indicate $A$ and $C$. The zero element ${}_A0_C=0\in\E(C,A)$ is called the {\it split $\E$-extension}. For any pair of $\E$-extensions ${}_A\del_C$ and ${}_{A'}\del{'}_{C'}$, let $\delta\oplus \delta'\in\mathbb{E}(C\oplus C', A\oplus A')$ be the
element corresponding to $(\delta,0,0,{\delta}{'})$ through the natural isomorphism $\mathbb{E}(C\oplus C', A\oplus A')\simeq\mathbb{E}(C, A)\oplus\mathbb{E}(C, A')
\oplus\mathbb{E}(C', A)\oplus\mathbb{E}(C', A')$.

For any $a\in\C(A,A')$ and $c\in\C(C',C)$,  $\E(C,a)(\del)\in\E(C,A')\ \ \text{and}\ \ \E(c,A)(\del)\in\E(C',A)$ are simply denoted by $a_{\ast}\del$ and $c^{\ast}\del$, respectively.

Let ${}_A\del_C$ and ${}_{A'}\del{'}_{C'}$ be any pair of $\E$-extensions. A {\it morphism} $(a,c)\colon\del\to{\delta}{'}$ of extensions is a pair of morphisms $a\in\C(A,A')$ and $c\in\C(C,C')$ in $\C$, satisfying the equality
$a_{\ast}\del=c^{\ast}{\delta}{'}$.}
Then the functoriality of $\E$ implies $\E(c,a)=a_{\ast}(c^{\ast}\del)=c^{\ast}(a_{\ast}\del)$.

\begin{definition}\cite[Definition 2.7]{HLN}
Let $\bf{C}_{\C}$ be the category of complexes in $\C$. As its full subcategory, define $\CC$ to be the category of complexes in $\C$ whose components are zero in the degrees outside of $\{0,1,\ldots,n+1\}$. Namely, an object in $\CC$ is a complex $X_{\bullet}=\{X_i,d_i^X\}$ of the form
\[ X_0\xrightarrow{d_0^X}X_1\xrightarrow{d_1^X}\cdots\xrightarrow{d_{n-1}^X}X_n\xrightarrow{d_n^X}X_{n+1}. \]
We write a morphism $f_{\bullet}\co X_{\bullet}\to Y_{\bullet}$ simply $f_{\bullet}=(f_0,f_1,\ldots,f_{n+1})$, only indicating the terms of degrees $0,\ldots,n+1$.
\end{definition}

\begin{definition}\cite[Definition 2.11]{HLN}
By Yoneda lemma, any extension $\del\in\E(C,A)$ induces natural transformations
\[ \del\ssh\colon\C(-,C)\Rightarrow\E(-,A)\ \ \text{and}\ \ \del\ush\colon\C(A,-)\Rightarrow\E(C,-). \]
For any $X\in\C$, these $(\del\ssh)_X$ and $\del\ush_X$ are given as follows.
\begin{enumerate}
\item[\rm(1)] $(\del\ssh)_X\colon\C(X,C)\to\E(X,A)\ :\ f\mapsto f\uas\del$.
\item[\rm (2)] $\del\ush_X\colon\C(A,X)\to\E(C,X)\ :\ g\mapsto g\sas\delta$.
\end{enumerate}
We simply denote $(\del\ssh)_X(f)$ and $\del\ush_X(g)$ by $\del\ssh(f)$ and $\del\ush(g)$, respectively.
\end{definition}

\begin{definition}\cite[Definition 2.9]{HLN}
 Let $\C,\E,n$ be as before. Define a category $\AE:=\AE^{n+2}_{(\C,\E)}$ as follows.
\begin{enumerate}
\item[\rm(1)]  A pair $\Xd$ is an object of the category $\AE$ with $X_{\bullet}\in\CC$
and $\del\in\E(X_{n+1},X_0)$, called an $\E$-attached
complex of length $n+2$, if it satisfies
$$(d_0^X)_{\ast}\del=0~~\textrm{and}~~(d^X_n)^{\ast}\del=0.$$
We also denote it by
$$X_0\xrightarrow{d_0^X}X_1\xrightarrow{d_1^X}\cdots\xrightarrow{d_{n-2}^X}X_{n-1}
\xrightarrow{d_{n-1}^X}X_n\xrightarrow{d_n^X}X_{n+1}\overset{\delta}{\dashrightarrow}.$$
\item[\rm (2)]  For such pairs $\Xd$ and $\langle Y_{\bullet},\rho\rangle$,  $f_{\bullet}\colon\Xd\to\langle Y_{\bullet},\rho\rangle$ is
defined to be a morphism in $\AE$ if it satisfies $(f_0)_{\ast}\del=(f_{n+1})^{\ast}\rho$.

\end{enumerate}
\end{definition}

\begin{definition}\cite[Definition 2.13]{HLN}\label{def1}
 An {\it $n$-exangle} is an object $\Xd$ in $\AE$ that satisfies the listed conditions.
\begin{enumerate}
\item[\rm (1)] The following sequence of functors $\C\op\to\Ab$ is exact.
$$
\C(-,X_0)\xrightarrow{\C(-,\ d^X_0)}\cdots\xrightarrow{\C(-,\ d^X_n)}\C(-,X_{n+1})\xrightarrow{~\del\ssh~}\E(-,X_0)
$$
\item[\rm (2)] The following sequence of functors $\C\to\Ab$ is exact.
$$
\C(X_{n+1},-)\xrightarrow{\C(d^X_n,\ -)}\cdots\xrightarrow{\C(d^X_0,\ -)}\C(X_0,-)\xrightarrow{~\del\ush~}\E(X_{n+1},-)
$$
\end{enumerate}
In particular any $n$-exangle is an object in $\AE$.
A {\it morphism of $n$-exangles} simply means a morphism in $\AE$. Thus $n$-exangles form a full subcategory of $\AE$.
\end{definition}

\begin{definition}\cite[Definition 2.22]{HLN}
Let $\s$ be a correspondence which associates a homotopic equivalence class $\s(\del)=[{}_A{X_{\bullet}}_C]$ to each extension $\del={}_A\del_C$. Such $\s$ is called a {\it realization} of $\E$ if it satisfies the following condition for any $\s(\del)=[X_{\bullet}]$ and any $\s(\rho)=[Y_{\bullet}]$.
\begin{itemize}
\item[{\rm (R0)}] For any morphism of extensions $(a,c)\co\del\to\rho$, there exists a morphism $f_{\bullet}\in\CC(X_{\bullet},Y_{\bullet})$ of the form $f_{\bullet}=(a,f_1,\ldots,f_n,c)$. Such $f_{\bullet}$ is called a {\it lift} of $(a,c)$.
\end{itemize}
In such a case, we simple say that \lq\lq$X_{\bullet}$ realizes $\del$" whenever they satisfy $\s(\del)=[X_{\bullet}]$.

Moreover, a realization $\s$ of $\E$ is said to be {\it exact} if it satisfies the following conditions.
\begin{itemize}
\item[{\rm (R1)}] For any $\s(\del)=[X_{\bullet}]$, the pair $\Xd$ is an $n$-exangle.
\item[{\rm (R2)}] For any $A\in\C$, the zero element ${}_A0_0=0\in\E(0,A)$ satisfies
\[ \s({}_A0_0)=[A\ov{\id_A}{\lra}A\to0\to\cdots\to0\to0]. \]
Dually, $\s({}_00_A)=[0\to0\to\cdots\to0\to A\ov{\id_A}{\lra}A]$ holds for any $A\in\C$.
\end{itemize}
Note that the above condition {\rm (R1)} does not depend on representatives of the class $[X_{\bullet}]$.
\end{definition}

\begin{definition}\cite[Definition 2.23]{HLN}
Let $\s$ be an exact realization of $\E$.
\begin{enumerate}
\item[\rm (1)] An $n$-exangle $\Xd$ is called an $\s$-{\it distinguished} $n$-exangle if it satisfies $\s(\del)=[X_{\bullet}]$. We often simply say {\it distinguished $n$-exangle} when $\s$ is clear from the context.
\item[\rm (2)]  An object $X_{\bullet}\in\CC$ is called an {\it $\s$-conflation} or simply a {\it conflation} if it realizes some extension $\del\in\E(X_{n+1},X_0)$.
\item[\rm (3)]  A morphism $f$ in $\C$ is called an {\it $\s$-inflation} or simply an {\it inflation} if it admits some conflation $X_{\bullet}\in\CC$ satisfying $d_0^X=f$.
\item[\rm (4)]  A morphism $g$ in $\C$ is called an {\it $\s$-deflation} or simply a {\it deflation} if it admits some conflation $X_{\bullet}\in\CC$ satisfying $d_n^X=g$.
\end{enumerate}
\end{definition}

\begin{definition}\cite[Definition 2.32]{HLN}
An {\it $n$-exangulated category} is a triplet $(\C,\E,\s)$ of additive category $\C$, additive bifunctor $\E\co\C\op\times\C\to\Ab$, and its exact realization $\s$, satisfying the following conditions.

(EA1) Let $A\ov{f}{\lra}B\ov{g}{\lra}C$ be any sequence of morphisms in $\C$. If both $f$ and $g$ are inflations, then so is $g\circ f$. Dually, if $f$ and $g$ are deflations, then so is $g\circ f$.

(EA2) For $\rho\in\E(D,A)$ and $c\in\C(C,D)$, let ${}_A\langle X_{\bullet},c\uas\rho\rangle_C$ and ${}_A\Yr_D$ be distinguished $n$-exangles. Then $(\id_A,c)$ has a {\it good lift} $f_{\bullet}$, in the sense that its mapping cone gives a distinguished $n$-exangle $\langle M^f_{\bullet},(d^X_0)\sas\rho\rangle$.

(EA2$\op$) Dual of {\rm (EA2)}.

Note that the case $n=1$, a triplet $(\C,\E,\s)$ is a  $1$-exangulated category if and only if it is an extriangulated category, see \cite[Proposition 4.3]{HLN}.
\end{definition}

\begin{example}
From \cite[Proposition 4.34]{HLN} and \cite[Proposition 4.5]{HLN},  we know that $n$-exact categories and $(n+2)$-angulated categories are $n$-exangulated categories.
There are some other examples of $n$-exangulated categories
 which are neither $n$-exact nor $(n+2)$-angulated, see \cite{HLN,HLN1,LZ,HZZ1}.
\end{example}

Let $(\C,\E,\s)$ be an $n$-exangulated category and  $\F\subseteq\E$ be an additive subfunctor {\rm (see \cite[Definition 3.7]{HLN})}. For a realization $\s$ of $\E$, define
$\s\hspace{-1.3mm}\mid_{\F}$ to be the restriction of $\s$ onto $\F$. Namely, it is defined by $\s\hspace{-1.3mm}\mid_{\F}(\del)=\s(\del)$ for any $\F$-extension $\del$.

\begin{lemma}\label{lem1}{\rm\cite[Proposition 3.16]{HLN}}
Let $(\C,\E,\s)$ be an $n$-exangulated category. For any additive subfunctor $\F\subseteq\E$, the following statements are equivalent.

{\rm (1)}~ $(\C,\F,\s\hspace{-1.2mm}\mid_{\F})$ is an $n$-exangulated category.

{\rm (2)}~ $\s\hspace{-1.2mm}\mid_{\F}$-inflations are closed under composition.

{\rm (3)}~ $\F\subseteq\E$ is closed.
\end{lemma}

\section{Hereditary $n$-exangulated categories}

We first recall the notion of \emph{cluster-tilting} subcategory from \cite{LZ1}.
\begin{definition}\label{defn}{\rm\cite[Definition 3.4]{LZ1}}
Let $\C$ be an $n$-exangulated category and $\X$ be an additive subcategory of $\C$.
$\X$ is called \emph{cluster-tilting} if

(1) $\E(\X,\X)=0$.

(2)  For any object $C\in\C$, there exist two distinguished $n$-exangles
$$X_0\xrightarrow{~~}X_1\xrightarrow{~~}\cdots\xrightarrow{~~}X_{n-1}\xrightarrow{~~}X_n\xrightarrow{~~}C\dashrightarrow$$
where $X_0, X_1,\cdots,X_{n}\in\X$ and
$$C\xrightarrow{~~}X'_1\xrightarrow{~~}X'_2\xrightarrow{~~}\cdots\xrightarrow{~~}X'_n\xrightarrow{~~}X'_{n+1}\dashrightarrow$$
where $X'_1,X'_2,\cdots,X'_{n+1}\in\X$.\end{definition}

\begin{remark}(1) An object $X$ is called cluster-tilting if $\rm add(X)$ is clsuter-tilting.

(2) When $\C$ is an $(n+2)$-angulated category, this definition is just \cite[Definition 5.3]{OT} and \cite[Definition 1.1]{ZZ1}.

\end{remark}

Now we give some examples of cluster-tilting subcategories.

\begin{example}\label{ex2}
{\upshape Let $\Lambda$ be the algebra given by the following (infinity) quiver with relations $x^2=0$:
	 \begin{align}
	 	\begin{minipage}{0.6\hsize}
	 		\ \ \ \ \  \xymatrix{\begin{smallmatrix}1\end{smallmatrix}&\begin{smallmatrix}2\end{smallmatrix}\ar[l]_{x}
	 			&\begin{smallmatrix}3\end{smallmatrix}\ar[l]_{x}&\begin{smallmatrix}4\end{smallmatrix}\ar[l]_{x}&\begin{smallmatrix}\cdots\end{smallmatrix}\ar[l]_{x}&\begin{smallmatrix}n\end{smallmatrix}\ar[l]_{x}&\begin{smallmatrix}\cdots\end{smallmatrix}\ar[l]_{x}&\begin{smallmatrix}\end{smallmatrix}}\notag
	 	\end{minipage}
	 \end{align}
The Auslander-Reiten quiver of ${\rm mod}\Lambda$ is the following:}
$$\xymatrix@C=0.2cm@R0.2cm{
&&&\bullet \ar[dr] &&\bullet \ar[dr]&&\bullet \ar[dr] &&\bullet \ar[dr] &&\bullet \ar[dr] &&\bullet \ar[dr] &&\bullet \ar[dr] &&\bullet \ar[dr]&&\bullet\ar[dr]\\
&&\bullet\ar[ur]  &&\circ\ar[ur]  &&\spadesuit\ar[ur]  &&\circ \ar[ur]  &&\bullet\ar[ur] &&\circ\ar[ur] &&\spadesuit \ar[ur]  &&\circ\ar[ur]   &&\bullet\ar[ur]   &&\cdots
}
$$
where the object denoted by $\spadesuit$ and $\bullet $ appear periodically. Let $\C$ be the additive closure of all the indecomposable objects denoted by $\spadesuit$ and $\bullet $.
 Then $\C$ is a cluster-tilting subcategory of ${\rm mod}\Lambda$,
hence it is $2$-abelian (see \cite[Theorem 3.16]{Ja}). Let $\X$ be the additive closure of all the indecomposable objects denoted by $\bullet$.
It is straightforward to verify that  $\X$ is a cluster-tilting subcategory of $\C$.
\end{example}

\begin{example}
Let $\Lambda$ be $n$-representation finite algebra, $\mathscr{O}_\Lambda$ the $(n+2)$-angulated cluster category associated to $\Lambda$. In particular, we let $n=3$ and $\mathscr{T}=\mathscr{O}_{A^3_2}$. This is the $5$-angulated (higher) cluster category of type $A_2$ (see \cite[Definition 5.2, Section 6 and Section 8]{OT}). The indecomposable objects can be identified with the elements of the set
$$\rm ind\mathscr{T}=\{1357, 1358, 1368, 1468, 2468, 2469, 2479, 2579, 3579 \}.$$
The Auslander-Reiten quiver of $\mathscr{T}$ is the following:
$$\begin{xy}
 (0,0)*+{\begin{smallmatrix}1368\end{smallmatrix}}="1",
(-15,-8)*+{\begin{smallmatrix}1468\end{smallmatrix}}="2",
(15,-8)*+{\begin{smallmatrix}1358\end{smallmatrix}}="3",
(-20,-20)*+{\begin{smallmatrix}2468\end{smallmatrix}}="4",
(20,-20)*+{\begin{smallmatrix}1357\end{smallmatrix}}="5",
(-15,-30)*+{\begin{smallmatrix} 2469\\
\end{smallmatrix}}="6",
(15,-30)*+{\begin{smallmatrix} 3579\\
\end{smallmatrix}}="7",
(-8,-42)*+{\begin{smallmatrix} 2479\\
\end{smallmatrix}}="8",
(8,-42)*+{\begin{smallmatrix} 2579\\
\end{smallmatrix}}="9",
\ar"1";"2", \ar"3";"1", \ar"2";"4", \ar"5";"3", \ar"4";"6",
\ar"7";"5",\ar"6";"8",\ar"9";"7",\ar"8";"9",
\end{xy}$$
It is straightforward to verify that the subcategory $$\mathcal T:={\rm add}(1357\oplus1358\oplus1368\oplus1468)$$ is cluster-tilting.
\end{example}

\begin{definition}\label{def2}\cite[Definition 3.14 and Definition 3.15]{ZW}
Let $(\C,\E,\s)$ be an $n$-exangulated category.
\begin{itemize}
\item[(1)] An object $P\in\C$ is called \emph{projective} if, for any distinguished $n$-exangle
$$A_0\xrightarrow{\alpha_0}A_1\xrightarrow{\alpha_1}A_2\xrightarrow{\alpha_2}\cdots\xrightarrow{\alpha_{n-2}}A_{n-1}
\xrightarrow{\alpha_{n-1}}A_n\xrightarrow{\alpha_n}A_{n+1}\overset{\delta}{\dashrightarrow}$$
and any morphism $c$ in $\C(P,A_{n+1})$, there exists a morphism $b\in\C(P,A_n)$ satisfying $\alpha_n\circ b=c$.
We denote the full subcategory of projective objects in $\C$ by $\P$.
Dually, the full subcategory of injective objects in $\C$ is denoted by $\I$.

\item[(2)] We say that $\C$ {\it has enough  projectives} if
for any object $C\in\C$, there exists a distinguished $n$-exangle
$$B\xrightarrow{\alpha_0}P_1\xrightarrow{\alpha_1}P_2\xrightarrow{\alpha_2}\cdots\xrightarrow{\alpha_{n-2}}P_{n-1}
\xrightarrow{\alpha_{n-1}}P_n\xrightarrow{\alpha_n}C\overset{\delta}{\dashrightarrow}$$
satisfying $P_1,P_2,\cdots,P_n\in\P$. We can define the notion of having \emph{enough injectives} dually.
\end{itemize}
\end{definition}
Next, we introduce the concept of a \emph{hereditary $n$-exangulated category}.
\begin{definition}\label{def1}
An $n$-exangulated category $(\C,\E,\s)$ is called \emph{hereditary} if for any object $C\in\C$, there exist a distinguished $n$-exangle
$$P_0\xrightarrow{~~}P_1\xrightarrow{~~}\cdots\xrightarrow{~~}P_{n-1}\xrightarrow{~~}P_n\xrightarrow{~~}C\dashrightarrow$$
where $P_0, P_1, \cdots, P_{n}\in\P$, or
$$C\xrightarrow{~~}I_1\xrightarrow{~~}I_2\xrightarrow{~~}\cdots\xrightarrow{~~}I_n\xrightarrow{~~}I_{n+1}\dashrightarrow$$
where $I_1,I_2,\cdots,I_{n+1}\in\I$.
\end{definition}

\begin{remark}
When $n=1$, the $1$-exangulated category $(\C,\E,\s)$  becomes an extriangulated category. If $(\C,\E,\s)$ has enough projectives and enough injectives, then the following statements are equivalent (see \cite[Propsition 2.1]{GNP}):

$\centerdot$ For any object $C\in\C$, there exists an $\E$-triangle $P_0\xrightarrow{~~}P_1\xrightarrow{~~}C\dashrightarrow$ with $P_0, P_1\in\P$.

$\centerdot$ For any object $C\in\C$, there exists an $\E$-triangle $C\xrightarrow{~~}I_1\xrightarrow{~~}I_2\dashrightarrow$ with $I_1, I_2\in\P$.

$\centerdot$ The additive bifunctor $\E^2(-,-)=0$.

Note that in \cite[Proposition 2.1]{GNP}, one of the key arguments in the proof is high-dimensional extension groups and long exact sequences. However, these may not necessarily exist in an $n$-exangulated category, so we just assume that one of the above two sequences exists in Definition \ref{def1}.
\end{remark}

Our first main result is as follows, which provides two examples of hereditary $n$-exangulated categories.

\begin{theorem}\label{main1}
Let $(\C,\E,\s)$ be an $n$-exangulated category and
$\X$ be a cluster-tilting subcategory of $\C$. We define
$$\hspace{4mm}\E_{\X}(C,A)=\{\delta\in\E(C,A)\mid(\del\ssh)_X=0~\mbox{for any}~X\in\X\},$$
$$\E^{\X}(C,A)=\{\delta\in\E(C,A)\mid\del\ush_X=0~\mbox{for any}~X\in\X\}.$$
Then

{\rm(1)} $(\C,\E_{\X},\s\hspace{-1.2mm}\mid_{\E_{\X}})$ is a hereditary $n$-exangulated category.

{\rm(2)} $(\C,\E^{\X},\s\hspace{-1.2mm}\mid_{\E^{\X}})$ is a hereditary $n$-exangulated category.
\end{theorem}
\proof
(1) By \cite[Proposition 3.19]{HLN}, we know that $\E_{\X}$ is a closed subfunctor of $\E$. Thus $(\A,\E_{\X},\s\hspace{-1.2mm}\mid_{\E_{\X}})$
is an $n$-exangulated category by Lemma \ref{lem1}. We only need to show that $(\A,\E_{\X},\s\hspace{-1.2mm}\mid_{\E_{\X}})$ is hereditary. Suppose that $$X_0\xrightarrow{}X_1\xrightarrow{}X_2\xrightarrow{}\cdots\xrightarrow{}X_{n-1}
\xrightarrow{}X_n\xrightarrow{g_n}X_{n+1}\overset{\delta}{\dashrightarrow}$$ is an any $\s\hspace{-1.2mm}\mid_{\E_{\X}}$-distinguished $n$-exangle.
For any object $Y_{n+1}\in\X$ and any morphism $f_{n+1}\in\C(Y_{n+1},X_{n+1})$, we have the following commutative diagram
$$\xymatrix{
X_0 \ar[r]^{}\ar@{=}[d]^{} & Y_1 \ar[r]^{}\ar[d]^{} & Y_2 \ar[r]^{} \ar[d]^{}& \cdots \ar[r]^{} & Y_{n} \ar[r]^{h_{n}}\ar[d]^{}&Y_{n+1}\ar[d]^{f_{n+1}}\ar@{-->}[r]^-{f_{n+1}\uas\delta} & \\
X_0 \ar[r]^{} & X_1 \ar[r]^{} & X_2\ar[r]^{} & \cdots \ar[r]^{} & X_n \ar[r]^{g_{n}} & X_{n+1}\ar@{-->}[r]^-{\del} & }$$
of $\s\hspace{-1.2mm}\mid_{\E_{\X}}$-distinguished $n$-exangles by $\rm (EA2)$. Since $\delta\in\E_{\X}(X_{n+1},X_0)$, we have $f_{n+1}\uas\delta=(\del\ssh)_X(f_{n+1})=0$. So there exists a morphism $h_{n+1}:Y_{n+1}\rightarrow X_n$, such that $f_{n+1}=g_{n}h_{n+1}$ by Lemma 3.3 in \cite{ZW}. This shows that every object in $\X$ is projective in $(\C,\E_{\X},\s\hspace{-1.2mm}\mid_{\E_{\X}})$. Since $\X$ is a cluster-tilting subcategory of $\C$, there exists an $\s$-distinguished $n$-exangle
$$ Z_0\xrightarrow{}Z_1\xrightarrow{}Z_2\xrightarrow{}\cdots\xrightarrow{}Z_{n-1}
\xrightarrow{}Z_n\xrightarrow{}C\overset{\varrho}{\dashrightarrow}\eqno{(\maltese)}$$
for any $C\in\C$, where $Z_0, Z_1, \cdots, Z_n\in \X$. For any $f:X\rightarrow C$ with $X\in\X$, we have the following commutative diagram
$$\xymatrix{
Z_0 \ar[r]^{}\ar@{=}[d]^{} & K_1 \ar[r]^{}\ar[d]^{} & K_2 \ar[r]^{} \ar[d]^{}& \cdots \ar[r]^{} & K_{n} \ar[r]^{}\ar[d]^{}&X\ar[d]^{f}\ar@{-->}[r]^-{f\uas\varrho} & \\
Z_0 \ar[r]^{} & Z_1 \ar[r]^{} & Z_2\ar[r]^{} & \cdots \ar[r]^{} & Z_n \ar[r]^{g_{n}} & C\ar@{-->}[r]^-{\varrho} & }$$
of $\s$-distinguished $n$-exangles by $\rm (EA2)$. Note that $Z_0, X\in\X$, then $f\uas\varrho=0$ since $\E(\X,\X)=0$. That is to say $(\varrho\ssh)_X(f)=f\uas\varrho=0$, so $\varrho\in \E_{\X}(C,Z_0)$. Thus the $\s\hspace{-1.2mm}\mid_{\E_{\X}}$-distinguished $n$-exangle $(\maltese)$ is also what we want in Definition \ref{def1}.

(2) It is similar to (1).
\qed
\vspace{1mm}

We recall the notion of $n$-extension closed subcategory from \cite{HLN}.

\begin{definition}\cite[Definition 2.34]{HLN}
Let $(\C,\E,\s)$ be an $n$-exangulated category. A subcategory $\A$ of $\C$ is called \emph{$n$-extension closed} if
 for any pair of objects $A$ and $C$ in $\A$ and any $\E$-extension $\delta\in\E(C,A)$, there is a distinguished $n$-exangle $\Xd$ with $X_i$ in $\A$ for $i=1,\cdots,n$.
\end{definition}

Let $(\C,\E,\s)$ be an $n$-exangulated category and $\A$ be an $n$-extension closed subcategory of $\C$. We define $\E_{\A}$ to be the restriction of
$\E$ onto $\A^{\rm op}\times\A$. For any $\delta\in\E_{\A}(C,A)$, take an $\s$-distinguished $n$-exangle $\Xd$ with $X_i$ in $\A$ for $i=1,\cdots,n$. We define $\t(\delta)=[X_{\bullet}]$, where the homotopy  equivalence class is taken in ${\bf{C}}^{n+2}_{({\A;\hspace{0.3mm}A,\hspace{0.3mm}C})}$.

\begin{proposition}\label{main1}
Let $(\C,\E,\s)$ be a hereditary $n$-exangulated category and
$\A$ be an $n$-extension closed subcategory of $\C$. Then

$(1)$  $(\A,\E_{\A},\t)$ is an $n$-exangulated category.

$(2)$ If $\P\subseteq\A$, where $\P$ denotes the full subcategory of projective objects in $(\C,\E,\s)$, then $(\A,\E_{\A},\t)$ is hereditary.
\end{proposition}

\proof $(1)$ By \cite[Theorem 3.3]{K}, we know that $(\A,\E_{\A},\t)$ is an $n$-exangulated category.

$(2)$ For any $P\in\P$, it is clear that $P$ is an projective object in $(\A,\E_{\A},\t)$. Let $A\in\A\subseteq\C$, since $(\C,\E,\s)$ is hereditary, there exists an $\s$-distinguished $n$-exangle
$$ P_0\xrightarrow{}P_1\xrightarrow{}P_2\xrightarrow{}\cdots\xrightarrow{}P_{n-1}
\xrightarrow{}P_n\xrightarrow{}A\overset{\delta}{\dashrightarrow}\eqno{(\diamondsuit)}$$
where $P_0,P_1,\cdots,P_n\in \P$.
This shows that $(\A,\E_{\A},\t)$ is hereditary.   \qed

We denote by $\C/\X$
the category whose objects are objects of $\C$ and whose morphisms are elements of
$\Hom_{\C}(A,B)/\X(A,B)$ for $A,B\in\C$, where $\X(A,B)$ is the subgroup of $\Hom_{\C}(A,B)$ consisting of morphisms
which factor through an object in $\X$.
$\C/\X$ is called the (additive) quotient category
of $\C$ by $\X$. For any morphism $f\colon A\to B$ in $\C$, we denote by $\overline{f}$ the image of $f$ under
the natural quotient functor $\C\to\C/\X$.
\medskip

Let $(\C, \mathbb{E}, \mathfrak{s})$ be an $n$-exangulated category and denote $\overline{\C}:=\C/\X$.
Assume that
$$A_0\xrightarrow{~\alpha_0~}A_1\xrightarrow{~\alpha_1~}A_2\xrightarrow{~\alpha_2~}\cdots\xrightarrow{~\alpha_{n-2}~}A_{n-1}
\xrightarrow{~\alpha_{n-1}~}A_n\xrightarrow{~\alpha_n~}A_{n+1}\overset{\delta}{\dashrightarrow}$$
is a distinguished $n$-exangle in $\C$. This sequence
$$A_0\xrightarrow{~\overline{\alpha_0}~}A_1\xrightarrow{~\overline{\alpha_1}~}A_2\xrightarrow{~\overline{\alpha_2}~}\cdots
\xrightarrow{~\overline{\alpha_{n-2}}~}A_{n-1}
\xrightarrow{~\overline{\alpha_{n-1}}~}A_n\xrightarrow{~\overline{\alpha_n}~}A_{n+1}$$
is called \emph{weak kernel-cokernel sequence} if the following sequences
$$
\overline{\C}(-,A_0)\xrightarrow{\overline{\C}(-,\  \overline{\alpha_0})}\overline{\C}(-,A_1)\xrightarrow{\overline{\C}(-,\ \overline{\alpha_1})}\cdots\xrightarrow{\overline{\C}(-,\ \overline{\alpha_{n-1}})}\overline{\C}(-,A_n)\xrightarrow{\overline{\C}(-,\ \overline{\alpha_n})}\overline{\C}(-,A_{n+1})
$$
and
$$
\overline{\C}(A_{n+1},-)\xrightarrow{\overline{\C}(\overline{{\alpha_n}},\ -)}\overline{\C}(A_{n},-)\xrightarrow{\overline{\C}(\overline{{\alpha_{n-1}}},\ -)}\cdots\xrightarrow{\overline{\C}(\overline{{\alpha_{1}}},\ -)}\overline{\C}(A_1,-)\xrightarrow{\overline{\C}(\overline{{\alpha_{0}}},\ -)}\overline{\C}(A_0,-).
$$
\begin{proposition}\label{main22}
Let $(\C, \mathbb{E}, \mathfrak{s})$ be a hereditary $n$-exangulated category and $\X$ a full subcategory of $\C$.
If $\X$ satisfies $\X\subseteq\P\cap\I$, then the ideal quotient $\C/\X$ is a hereditary $n$-exangulated category
if and only if any distinguished $n$-exangle in $\C$ induces a weak kernel-cokernel sequence in $\C/\X$.
\end{proposition}

\proof Define the additive
bifunctor
$\overline{\E}\colon \overline{\C}^{\rm op}\times\overline{\C}\to\Ab$
given by
\begin{itemize}
\item $\overline{\E}(C,A)=\E(C,A)$ for any $A,C\in\C$,

\item $\overline{\E}(\overline{c},\overline{a})=\E(c,a)$
for any $a\in\C(A,A'),~c\in\C(C,C')$, where $\overline{a}$ and $\overline{c}$
denote the images of $a$ and $c$ in $\C/\X$, respectively.
\end{itemize}

For any $\overline{\E}$-extension $\delta\in\overline{\E}(C,A)={ \E}(C,A)$, define
$$\overline{\s}(\delta)=\overline{\s(\delta)}=[A\xrightarrow{~\overline{\alpha_0}~}
B_1\xrightarrow{~\overline{\alpha_1}~}
B_2\xrightarrow{~\overline{\alpha_2}~}\cdots\xrightarrow{~\overline{\alpha_{n-1}}~}B_n\xrightarrow{~\overline{\alpha_{n}}~}C]$$
using $\s(\delta)=[A\xrightarrow{~\alpha_0~}
B_1\xrightarrow{~\alpha_1~}
B_2\xrightarrow{~\alpha_2~}\cdots\xrightarrow{~\alpha_{n-1}~}B_n\xrightarrow{\alpha_{n}}C]$.

By \cite[Theorem 3.1]{HZZ1}, we know that  $(\C/\X, \overline{\E}, \overline{\s})$ is an $n$-exangulated category
if and only if any distinguished $n$-exangle in $\C$ induces a weak kernel-cokernel sequence in $\C/\X$. We need to show that $\C/\X$ is hereditary. Let $C\in\C/\X$, since $(\C,\E,\s)$ is hereditary, there exists an $\s$-distinguished $n$-exangle
$$ P_0\xrightarrow{\alpha_0}P_1\xrightarrow{\alpha_1}P_2\xrightarrow{\alpha_2}\cdots\xrightarrow{}P_{n-1}
\xrightarrow{\alpha_{n-1}}P_n\xrightarrow{\alpha_{n}}C\overset{\delta}{\dashrightarrow},$$
where $P_0,P_1,\cdots,P_n\in \P$.
Thus $$ P_0\xrightarrow{\overline{\alpha_0}}P_1\xrightarrow{\overline{\alpha_1}}P_2\xrightarrow{\overline{\alpha_2}}\cdots\xrightarrow{}P_{n-1}
\xrightarrow{\overline{\alpha_{n-1}}}P_n\xrightarrow{\overline{\alpha_{n}}}C\overset{\overline{\delta}}{\dashrightarrow}\eqno{(\maltese\maltese)}$$
is an $\overline{\s}$-distinguished $n$-exangle. It is clear that $P_i$ are also projective objects in $\C/\X$, where $i=0,1,2,\cdots,n$. In fact, let $$A_0\xrightarrow{~\overline{\alpha_0}~}A_1\xrightarrow{~\overline{\alpha_1}~}A_2\xrightarrow{~\overline{\alpha_2}~}\cdots
\xrightarrow{~\overline{\alpha_{n-2}}~}A_{n-1}
\xrightarrow{~\overline{\alpha_{n-1}}~}A_n\xrightarrow{~\overline{\alpha_n}~}A_{n+1}\overset{\overline{\delta}}{\dashrightarrow}$$
be any $\overline{\s}$-distinguished $n$-exangle, we know that $$A_0\xrightarrow{~\alpha_0~}A_1\xrightarrow{~\alpha_1~}A_2\xrightarrow{~\alpha_2~}\cdots\xrightarrow{~\alpha_{n-2}~}A_{n-1}
\xrightarrow{~\alpha_{n-1}~}A_n\xrightarrow{~\alpha_n~}A_{n+1}\overset{\delta}{\dashrightarrow}$$
is an $\s$-distinguished $n$-exangle. For any morphism $\overline{f_i}\in\overline{\C}(P_i,A_{n+1})$, since $P_i$ are projective, there exists $g_i\in{\C}(P_i,A_{n})$, such that $\alpha_ng=f_i$ and then $\overline{\alpha_n}\circ\overline{g}=\overline{f_i}$. Thus the  $\overline{\s}$-distinguished $n$-exangle $(\maltese\maltese)$ is also what we want in Definition \ref{def1}.

 \qed

\section{$0$-Auslander $n$-exangulated categories}

We first introduce the concept of a \emph{$0$-Auslander $n$-exangulated category.}

\begin{definition}\label{def2}
An $n$-exangulated category $(\C,\E,\s)$ is called \emph{$0$-Auslander} if

(1) $\C$ is hereditary.

(2)  For any projective object $P\in\P$, there exists a distinguished $n$-exangle
$$P\xrightarrow{~~}Q_1\xrightarrow{~~}Q_2\xrightarrow{~~}\cdots\xrightarrow{~~}Q_{n-1}\xrightarrow{~~}Q_n\xrightarrow{~~}I\dashrightarrow$$
where $Q_1, Q_2, \cdots, Q_{n}$ are projective-injective objects and $I\in\I$.
\end{definition}

\begin{remark}
When $n=1$, the Definition \ref{def2} is slightly different from that in \cite[Definition 1.1]{GNP}, we do not assume that the $\C$ has enough projectives.
\end{remark}

Let $(\C,\Sigma^n,\Theta)$ be an $(n+2)$-angulated category. Since $\Sigma^n\colon\C\xrightarrow{~\simeq~}\C$ is an  automorphism, then
$\Sigma^n$ gives an additive bifunctor
$$\E_{\Sigma^n}=\C(-,\Sigma^n-)\colon \C^{\rm op}\times \C\to {\rm Ab},$$
defined by the following.
\begin{itemize}
\item[\rm (i)] For any $A,C\in\C$, $\E_{\Sigma^n}(C, A)=\C(C,\Sigma A)$;

\item[\rm (ii)] For any $a\in\C(A,A')$ and $c\in\C(C',C)$, the map $\E_{\Sigma^n}(c, a)\colon\C(C, \Sigma^n A)\to \C(C', \Sigma^n A)$
sends $\delta\in\C(C, \Sigma^n A)$ to $c^{\ast}a_{\ast}\delta=(\Sigma^n a)\circ\delta\circ c$.
\end{itemize}

For each $\delta\in\E_{\Sigma^n}(C, A)$, we complete
it into an $(n+2)$-angle
$$A\xrightarrow{f_0}X_1\xrightarrow{f_1}X_2\xrightarrow{f_2}\cdots\xrightarrow{f_{n-1}}X_n\xrightarrow{f_n}C\xrightarrow{\delta}\Sigma^n A_0$$
Define $\s_{\Theta}(\delta)=[X_{\bullet}]$ by using $X_{\bullet}\in{\mathbf{C}^{n+2}_{(A,\hspace{0.8mm}C)}}$ given by
$$A\xrightarrow{f_0}X_1\xrightarrow{f_1}X_2\xrightarrow{f_2}\cdots\xrightarrow{f_{n-1}}X_n\xrightarrow{f_n}C$$

The following result shows that any $(n+2)$-angulated category can be viewed as
an $n$-exangulated category.
\begin{theorem}{\rm \cite[Proposition 4.5]{HLN}}
With the above definition, $(\C,\E_{\Sigma^n},\s_{\Theta})$ is an $n$-exangulated category.
\end{theorem}
Keep the mark in front, let $\T$ be a cluster-tilting subcategory of $(\C,\Sigma^n,\Theta)$. We define
$$\F_{\T}(C,A)=\{A\xrightarrow{~~}X_1\xrightarrow{~~}\cdots\xrightarrow{~~}
X_{n-1}\xrightarrow{~~}X_n\xrightarrow{~~}C\overset{\delta}{\dashrightarrow}~\mid~\del~ \text{is induced by the}~(n+2)\text{-angle} $$
$$A\xrightarrow{~~}X_1\xrightarrow{~~}\cdots\xrightarrow{~~}X_{n-1}
\xrightarrow{~~}X_n\xrightarrow{~~}C\xrightarrow{~f~}\Sigma^nA~\text{with}~f~\text{factors through}~\Sigma^nT,\hspace{0.3mm}T\in\T\}. $$

Let $\C$ be  an additive category and $\X$ be a subcategory of $\C$.
Recall that we say a morphism $f\colon A \to B$ in $\C$ is an $\X$-\emph{monic} if
$${\C}(f,X)\colon {\C}(B,X) \to {\C}(A,X)$$
is an epimorphism for all $X\in\X$. We say that $f$ is an $\X$-\emph{epic} if
$${\C}(X,f)\colon {\C}(X,A) \to {\C}(X,B)$$
is an epimorphism for all $X\in\X$.

Our second main result is as follows,
which provides an example of $0$-Auslander $n$-exangulated categories.

\begin{theorem}\label{main2}
Let $(\C,\Sigma^n,\Theta)$ be an $(n+2)$-angulated category and
$\T$ be a cluster-tilting subcategory of $\C$.
Then

$(1)$ $(\C,\F_{\T},\s_{\Theta}\hspace{-1.2mm}\mid_{\F_{\T}})$
is an $n$-exangulated category.

$(2)$ If $\P\subseteq\T$, where $\P$ denotes the full subcategory of projective objects in $(\C,\F_{\T},\s_{\Theta}\hspace{-1.2mm}\mid_{\F_{\T}})$, then $(\C,\F_{\T},\s_{\Theta}\hspace{-1.2mm}\mid_{\F_{\T}})$
is $0$-Auslander.

\end{theorem}
\proof
$(1)$ Let $$A\xrightarrow{~f_0~}X_1\xrightarrow{~~}X_2\xrightarrow{~~}\cdots\xrightarrow{~~}X_{n-1}\xrightarrow{~~}X_n\xrightarrow{~~}C\overset{\delta}\dashrightarrow$$
in $\F_{\T}(C,A)$. Then we have an $(n+2)$-angle
$$A\xrightarrow{~f_{0}~}X_1\xrightarrow{~~}\cdots\xrightarrow{~~}X_{n-1}
\xrightarrow{~~}X_n\xrightarrow{~~}C\xrightarrow{~f~}\Sigma^n A$$
with $f$ factors through $\Sigma^n T$, where $T\in\T$. That is, there exist two morphisms $g\colon C\rightarrow \Sigma^n T$ and $h\colon\Sigma^n T\rightarrow \Sigma^n A$ such that $f=hg$. We claim that ${f_0}$ is $\Sigma^n \T$-monic. Indeed, we can obtain the following $(n+2)$-angle
$\Sigma^{-n} C\xrightarrow{}A\xrightarrow{~f_{0}~}X_1\xrightarrow{}\cdots\xrightarrow{~~}
X_{n-1}\xrightarrow{~~}X_n\xrightarrow{~~}C.$

For any $T'\in\T$ and $k:A\rightarrow \Sigma^n T'$, we have the following commutative diagram
$$\xymatrix{
\Sigma^{-n} C\ar[r]^{\Sigma^{-n}f}\ar[d]_{\Sigma^{-n}g} & A \ar[r]^{f_{0}}\ar[d]^{k} & X_1 \ar@{-->}[dl]^{m} \ar[r]^{} & \cdots \ar[r]^{} & X_{n-1} \ar[r]^{}&X_{n}\ar[r]^{} &C. &\\
T \ar[ur]^{\Sigma^{-n}h} & \Sigma^{n}T'  &  &   &  &  && }$$
Since $\T$ be a cluster-tilting subcategory of $\C$, we have $\C(\T,\Sigma^n\T)=0$
which implies $$k\circ\Sigma^{-n}f=k\circ\Sigma^{-n}h\circ\Sigma^{-n}g=0\circ\Sigma^{-n}g=0.$$ There exists a morphism $m:X_1 \rightarrow \Sigma^n T'$, such that $k=mf_0$. Therefore, $f_0$ is $\Sigma^n \T$-monic. By \cite[Lemma 4.2]{LZ}, we have that $\F_{\T}$ is an additive subfunctor of $\E_{\Sigma^n}$ and $\Sigma^n \T$-monomorphisms are closed under composition. Hence $(\C,\F_{\T},\s_{\Theta}\hspace{-1.2mm}\mid_{\F_{\T}})$
is an $n$-exangulated category by Lemma \ref{lem1}.

(2) In order to prove $(\C,\F_{\T},\s_{\Theta}\hspace{-1.2mm}\mid_{\F_{\T}})$ is $0$-Auslander, we first prove that $\Sigma^n T$ is injective object and $T$ is projective object in  $(\C,\F_{\T},\s_{\Theta}\hspace{-1.2mm}\mid_{\F_{\T}})$ for any $T\in\T$. In fact, the above proof also shows that $\Sigma^n T$ is a projective object. Let
$$A\xrightarrow{}X_1\xrightarrow{~~}X_2\xrightarrow{~~}\cdots\xrightarrow{~~}X_{n-1}
\xrightarrow{~~}X_n\xrightarrow{~g~}C\overset{\delta}\dashrightarrow$$
be an $\s_{\Theta}\hspace{-1.2mm}\mid_{\F_{\T}}$-distinguished $n$-exangle. Then we have an $(n+2)$-angle
$$A\xrightarrow{}X_1\xrightarrow{~~}\cdots\xrightarrow{~~}X_{n-1}\xrightarrow{~~}
X_n\xrightarrow{~g~}C\xrightarrow{~h~}\Sigma^n A,$$
where $h$ factors through $\Sigma^n T$ with $T\in\T$. For any morphism $f:T'\rightarrow C$, where $T'\in\T$, we have the following commutative diagram
$$\xymatrix{&   &  &   &  & T'\ar[d]^{f}\ar@{-->}[dl]_{k} &\Sigma^n T\ar[d]^{q}& \\
A\ar[r]^{} & X_1 \ar[r]^{} & X_2 \ar[r]^{} & \cdots \ar[r]^{} & X_{n} \ar[r]^{g}&C\ar[r]^{h} \ar[ur]^{m}&\Sigma^nA. &}$$
Note that $\C(\T,\Sigma^n\T)=0$, so $hf=qmf=q0=0$. There exists a morphism $k\colon T' \rightarrow X_n$ such that $f=gk$. Therefore, $T'$ is projective.
Since $\T$ be a cluster-tilting subcategory of $\C$, there exists an $(n+2)$-angle
$$ T_0\xrightarrow{}T_1\xrightarrow{}T_2\xrightarrow{}\cdots\xrightarrow{}T_{n-1}
\xrightarrow{}T_n\xrightarrow{}C\xrightarrow{~f~}\Sigma^nT_0$$
for any $C\in\C$, where $T_0,T_1,\cdots,T_n\in\T$.
Since we have the following commutative diagram
$$\xymatrix{&   &  &   &  &  &\Sigma^n T_0\ar@{=}[d]^{}& \\
T_0\ar[r]^{} & T_1 \ar[r]^{} & T_2 \ar[r]^{} & \cdots \ar[r]^{} & T_{n} \ar[r]^{}&C\ar[r]^{f} \ar[ur]^{f}&\Sigma^nT_0. &}$$
So we obtain an $\s_{\Theta}\hspace{-1.2mm}\mid_{\F_{\T}}$-distinguished $n$-exangle $$T_0\xrightarrow{}T_1\xrightarrow{~~}T_2\xrightarrow{~~}\cdots\xrightarrow{~~}T_{n-1}\xrightarrow{~~}T_n\xrightarrow{g}C\overset{\delta}\dashrightarrow.\eqno{(\heartsuit)}$$
The $\s\hspace{-1.2mm}\mid_{\E_{\X}}$-distinguished $n$-exangle $(\heartsuit)$ is what we want in Definition \ref{def1}. For any $P\in\P$, there exists an $(n+2)$-angle
$$P\xrightarrow{}0\xrightarrow{~~}\cdots\xrightarrow{~~}0\xrightarrow{~~}0\xrightarrow{}\Sigma^n P\xrightarrow{}\Sigma^n P.$$
Note that $\T=\P$, then $P\in\T$. We have the following commutative diagram
$$\xymatrix{&   &  &   &  &  &\Sigma^n P\ar[d]^{}& \\
P\ar[r]^{} & 0 \ar[r]^{} & 0 \ar[r]^{} & \cdots \ar[r]^{} & 0 \ar[r]^{}&\Sigma^nP\ar[r]^{} \ar[ur]^{}&\Sigma^nP.}$$
So we obtain an $\s_{\Theta}\hspace{-1.2mm}\mid_{\F_{\T}}$-distinguished $n$-exangle $$P\xrightarrow{}0\xrightarrow{~~}0\xrightarrow{~~}\cdots\xrightarrow{~~}0
\xrightarrow{~~}0\xrightarrow{~g~}\Sigma^nP\overset{}\dashrightarrow.\eqno{(\heartsuit\heartsuit)}$$
Since the object $0$ is projective-injective objects and $\Sigma^nP\in\I$, hence the $\s\hspace{-1.2mm}\mid_{\E_{\X}}$-distinguished $n$-exangle $(\heartsuit\heartsuit)$ is what we want in Definition \ref{def2}(2). This shows that $(\C,\F_{\T},\s_{\Theta}\hspace{-1.2mm}\mid_{\F_{\T}})$
is $0$-Auslander.  \qed

\textbf{Jian He}\\
Department of Applied Mathematics, Lanzhou University of Technology, 730050 Lanzhou, Gansu, P. R. China\\
E-mail: \textsf{jianhe30@163.com}\\[0.3cm]
\textbf{Jing He}\\
School of Science, Hunan University of Technology and Business, 410205 Changsha, Hunan, P. R. China\\
E-mail: \textsf{jinghe1003@163.com}\\[0.3cm]
\textbf{Panyue Zhou}\\
School of Mathematics and Statistics, Changsha University of Science and Technology, 410114 Changsha, Hunan,  P. R. China\\
E-mail: \textsf{panyuezhou@163.com}

\end{document}